\documentclass[9pt,amstex]{article}

\usepackage{amssymb,amsmath} 
\usepackage{latexsym}

\newtheorem{thm}{Theorem}[section] 

\newtheorem{rmk}[thm]{Remark}

\newtheorem{lemma}[thm]{Lemma}

\def\cyclic{\mathop{\kern0.9ex{{+}
\kern-2.2ex\raise-.28ex\hbox{\Large\hbox
{$\circlearrowright$}}}}\limits}

\def\buildrel#1_#2^#3{\mathrel{\mathop{\kern 0pt#1}\limits_{#2}^{#3}}}

\newcommand{\Pf}{\noindent{\em Proof.} }
\newcommand{\EPf}
{%
\mbox{}%
\nolinebreak%
\hfill%
\rule{2mm}{2mm}%
\medbreak%
\par%
}

\newcommand{\im}{\mbox{$\mathtt{Im}$}}

\newcommand{\C}{\mathbb C} 

\newcommand{\A}{\mathbb A} 
\newcommand{\B}{\mathbb B}

\newcommand{\R}{\mathbb R}
\newcommand{\K}{\mathbb K}
 
\renewcommand{\H}{\mathbb H} 
\newcommand{\N}{\mathbb N}

\newcommand{\g}{{\mathfrak{g}}{}}

\newcommand{\Der}{{\mathfrak{Der}}{}} 
\newcommand{\q}{{\mathfrak{q}}{}}

\newcommand{\h}{{\mathfrak{h}}{}} 
\renewcommand{\sp}{{\mathfrak{sp}}{}}

\def\cref#1{Corollary~\ref{#1}}

\newcommand{\sP}{\mathfrak{sp}}
\newcommand{\oo}{\mathfrak{o}}
\newcommand{\so}{\mathfrak{so}}
\newcommand{\gL}{\mathfrak{gl}}
\newcommand{\uu}{\mathfrak{u}}

\newcommand{\m}{\mathfrak{m}}

\newcommand{\Gl}{\mathrm{Gl}}
\newcommand{\UU}{\mathrm{U}}
\newcommand{\OO}{\mathrm{O}}

\newcommand{\Sp}{\mathrm{Sp}}
\newcommand{\Herm}{\mathrm{Herm}}
\newcommand{\Sym}{\mathrm{Sym}}
\newcommand{\Aherm}{\mathrm{Aherm}}
\newcommand{\Asym}{\mathrm{Asym}}

\newcommand{\Gras}{\mathrm{Gras}}

\newcommand{\id}{\mathrm{id}}

\newcommand{\cX}{\mathcal{X}}

\newcommand{\bs}{\mathbf{s}}
\newcommand{\bt}{\mathbf{t}}

\newcommand{\HHH}{\mathbb H}

\newcommand{\bF}{\mathbb{F}}
\newcommand{\msk}{\medskip}
\newcommand{\ssk}{\smallskip}
\newcommand{\nin}{\noindent}

\title{Homotopes of Symmetric Spaces\\
I. Construction by Algebras with Two Involutions}

\author{
{\bf Wolfgang Bertram}\\
IECN Universit\'e Nancy I, France\\
e-mail:  {\tt bertram@iecn.u-nancy.fr}\\
{\bf Pierre Bieliavsky}\\
Universit\'e Catholique de Louvain, Belgium.\\
e-mail: {\tt bieliavsky@math.ucl.ac.be}\\
}
\date{}

\begin{document}
\maketitle

{\small
\nin {\bf Abstract.}
We investigate a special kind of contraction of symmetric spaces (respectively, of Lie
triple systems), called {\em homotopy}.
In this first part of a series of two papers we construct such contractions for 
classical symmetric spaces in an elementary way by using {\em associative algebras
with several involutions}. This construction shows a remarkable duality between
the underlying ``space'' and the ``deformation parameter''.

\msk \nin
 
\msk \nin
}

\section*{Introduction}

Usually, ``classical groups'' and ``classical symmetric spaces'' come in 
discrete families, parametrized by natural numbers $n \in \N$ (see \cite{Berger57}
for the classification of simple symmetric spaces).
 In this work, which is the first
in a series of two papers (mainly independent of each other), we will show that
there are indeed ``continuous families'' of such spaces. For classical groups,
such as $\Gl_n(\K)$, this has been observed in \cite{BeKi09}: 
fix $n \in \N$, then there is a ``continuous family of groups'', parametrized
by matrices $A \in M(n,n;\K)$,
$$
G_A : =\{ X \in M(n,n;\K) \mid \det (1 - XA) \not= 0 \}, 
\quad \mbox{ with product } \quad
X \cdot_A Y := X + Y - XAY .
$$
For $A=1$, the unit matrix, this is just $\Gl_n(\K)$ (with origin translated to zero),
and for $A=0$, it is the abelian vector group of matrices. 
In this first part of our work  we will give such direct linear algebra constructions
for certain contractions, called {\em homotopes}, of all ``classical symmetric spaces''.
In the following second part, using Jordan algebra techniques, we will show
that the lists thus obtained are indeed {\em complete} lists of such homotopes. 
In both papers, we take  an algebraic viewpoint, working with
{\em Lie triple systems} (which are infinitesimal versions of symmetric spaces).
We intend to use these constructions in further work to generalize quantization procedures
from \cite{BDS}.

\ssk
The basic idea of a direct construction of homotopes of classical Lie triple systems
is as follows: in the setting of an associative algebra
$\A$, ``homotopy'' amounts to the observation that,
for any $A \in \A$, the product $(X,Y) \mapsto XAY$ on $\A$ is again associative, hence
$[X,Y]_A := XAY - YAX$ is a Lie bracket, and 
\begin{equation}
\label{0.2}
[X,Y,Z]_A:= [[X,Y]_A,Z]_A =(XAYAZ+ZAYAX) -(YAXAZ+ZAXAY)
\end{equation}
is a Lie triple product on $\A$. For $A=1$ (unit element), 
this is the ``standard'' or ``general linear''
Lie triple product, and for
$A=0$ it is the ``flat'' Lie triple product on $\A$; thus we may say that the family 
of Lie triple products indexed by $A \in \A$ is a contraction of the ``general linear'' Lie triple
product. 
Now assume that the
 associative algebra $\A$ carries several commuting involutions $\tau_1,\ldots,
\tau_k$; for most of the constructions
 $k=1$ or $2$ will be sufficient. Let $\m$ be any of the $2^k$ joint eigenspaces
of these involutions. The key observation is now: if $A$  belongs to some joint eigenspace, then
$\m$ is stable under $[X,Y,Z]_A$, hence is a Lie triple system (Lemma \ref{TwoInv}). 
This means that we have $4^k$ families of Lie triple systems, which can be organized in a
table forming  a ``$2^k \times 2^k$-matrix'', with rows containing homotopes living on the same
underlying space $\m$ and columns containing homotopes parametrized by the same
parameter space. Such kind of ``duality'' between space and deformation parameter is
a special feature of the classical spaces governed by {\em associative} algebras considered here;
it is not present in the general (Jordan theoretic) formulation to be discussed in Part II of this work.

\ssk
We give several examples of such ``matrix-tables'' (Theorems \ref{ProjHomotopeTh},
\ref{SiegelTh}, \ref{QuaternionicTh1}, \ref{QuaternionicTh2}); all of them
reflect interesting features of families of classical symmetric spaces, among them, e.g.,
the {\em Siegel upper half plane} $\Sp_n(\R)/\UU(n)$ (Theorem \ref{SiegelTh}) and its
compact dual $\Sp(n)/\UU(n)$ (Theorem \ref{QuaternionicTh1}). 
The  duality between lines and columns of these tables makes it now easy 
 to calculate an {\em algebra imbedding} of
the Lie triple system $\m$, that is, a Lie algebra $\g$ such that $\g =
\h \oplus \m$ is a symmetric Lie algebra decomposition: 
if $\m$ belongs to the antidiagonal of the table, then $\m$ is already a Lie algebra;
else we may choose for $\h$ the space from the antidiagonal that belongs to the same
column as $\m$ (see Theorem \ref{TwoInvTh} and its proof).

\ssk
For a fixed joint eigenspace $\m$,
the Lie triple systems $\m$ with Lie triple product  $[X,Y,Z]_A$ can be interpreted as
curvature tensors of symmetric spaces $M =G/H$ depending on $A$ and all having the
same tangent space at the origin.
 In Section 4, we focus on the real, finite-dimensional case and give a list of all families of symmetric
spaces obtained that way, organized according to the point of view of 
``deformations'' and ``contractions'':
for each underlying space $\m$ we give a list of Lie triple products on $\m$ that are homotopes
of each other (Theorem \ref{ClassificationTheorem}).

\ssk
This work is organized as follows:
basic facts on Lie triple systems and classical Lie algebras are recalled in Chapters 1 and 2;
in Chapter 3 the ``two-involution-construction'' is explained, and the most important 
examples are woked out. Chapter 4 contains the systematic list  of contractions
of symmetric spaces thus obtained.
In Part II of this work (\cite{BeBi}) we explain a more general construction of 
homotopes, which uses basic ideas of {\em Jordan theory} (as we hope to convince
the reader, the concept of homotopy is useful in the associative theory, but its r\^ole in
the non-associative theory is certainly even more important): we define the
{\em structure variety}, which is the natural parameter space for contractions, and
we show that the list given here in Chapter 4 is, essentially, a complete description
of the corresponding structure varieties. While proving this, we will also obtain results
describing in more detail the structure of the contracted spaces: a
non-reductive homotope has a {\em bundle structure}, with flat fibers and 
a reductive base, and we will analyze some interesting low-dimensional examples of such
fibered symmetric spaces.

\msk
{\bf Acknowledgements.}
W.B.\ thanks Universit\'e Catholique de Louvain
 for hospitality in 2010 when part of this work was carried out.
P.B. thanks  Universit\'e Henri Poincar\'e-Nancy I for hospitality and
 the Belgian Scientific Policy (BELSPO) for its support through the 
IAP `NOSY' to which he is affiliated at the Universit\'e Catholique de Louvain.

\ssk
{\bf Notation.} Throughout this paper $\K$ is a commutative base ring in which $2$ is invertible.

\section{Lie triple systems} \label{sec:Liealg}

A {\em Lie triple system (LTS)} is a $\K$-module $\q$ together with a trilinear map
$$
\q^3 \to \q, \quad (X,Y,Z) \mapsto [X,Y,Z] =:R(X,Y)Z
$$
satisfying, for all $X,Y,Z,U,V \in \q$,

\begin{description}
\item[(LT1)] $[X,Y,Z] = - [Y,X,Z]$
\item[(LT2)] $[X,Y,Z]+[Y,Z,X] + [Z,X,Y]=0$
\item[(LT3)] the endomorphism $D:=R(U,V)$ is a derivation of the trilinear product
$[X,Y,Z]$. 
\end{description}

\nin
Every Lie algebra $\g$ with $[X,Y,Z]:=[[X,Y],Z]$ is a LTS, and if $\sigma$ is an automorphism
of $\g$ of order $2$, then
the $-1$-eigenspace $\q$ of $\sigma$ is stable under this trilinear product and hence is a LTS.
Every LTS $\q$ is obtained in this way: we may take for $\g$ the {\em standard imbedding}
$\q \oplus [\q,\q] \subset \q \oplus \Der(\q)$ (see \cite{Lo69}).
The pair $(\g,\sigma)$ is called a {\em symmetric pair}. 
If $\h = \g^\sigma$ is the fixed point algebra of $\sigma$, we will sometimes, with some abuse
of notation, denote the symmetric pair also by $(\g,\h)$. This notation is motivated (in the real
finite dimensional case) by the usual
description of the {\em associated symmetric  space} as a homogeneous space
$M=G/H$, where $G$ is a Lie group with involution $\sigma$
and $H$ an open subgroup of the fixed point group $G^\sigma$.
The case of a Lie group $H$ considered as symmetric space
$H \times H/diag(H \times H)$ will be called a {\em group case}; it corresponds to the case
of a Lie algebra, considered as a LTS.

\paragraph{Cartan-duality.}
It is clear that, if $(\q,R)$ is a LTS, then all multiples $(\q,\lambda R)$ for $\lambda \in \K$
are again LTS. In particular,
$(\q,-R)$ is again a LTS, called the
{\em $c$-dual Lie triple system}, where the letter $c$ refers to ``compact'' or ``Cartan'':
indeed, in the real finite dimensional case, $R$ is of compact type if and only if
$-R$ is of non-compact type. Note that, if $\lambda$ is a square in $\K^\times$,
then $\lambda \id_\q$ is an isomorphism between $R$ and $\lambda^2 R$.
In particular, for $\K=\C$, the LTS $R$ and  $-R$ are always isomorphic to each other.
For $\K=\R$, if $\g = \h \oplus \q$ is the standard imbedding of $R$, then
$\h \oplus i \q$ (subalgebra of $\g_\C$) is the standard imbedding of $-R$.
Using this, it is seen, for example, that the $c$-dual of a group case $H$ is a symmetric
space of the form $H_\C/H$, where $H_\C$ is a complexification of $H$.

\section{Classical Lie algebras}\label{ClassLieAlgebras}

We call ``classical'' Lie algebras that are defined by means of
 an associative unital $\K$-algebra $\A$  with some
involution  $\tau$ (antiautomorphism of order $2$); we often use the notation $a^*$ for $\tau(a)$
and  denote
the eigenspace decomposition of $*$ by
$\A = \A^\tau \oplus \A^{-\tau}= \Herm(\A,*) \oplus \Aherm(\A,*)$. 
We are going to define families of classical Lie algebras associated to these data, 
parametrized by certain elements $A \in \A$.
The three main types of classical Lie algebras are given by

\begin{lemma} \label{algLemma}
The following data define Lie algebras:
\begin{description}
\item[(1)]   {\em general linear:}
$\g_A:=\g_{\A,A}:=\A$ with 
$[X,Y]_A:=XAY-YAX$, for any $A \in \A$,
\item[(2)]  {\em unitary / orthogonal:}
$\uu_A := \uu_{\A,A,\tau}:=\Aherm(\A,*)$ with 
$[X,Y]_A$, for any $A \in \Herm(A,*)$,
\item[(3)]
{\em (half) symplectic:}
$\sp_A := \sp_{\A,A,\tau}:=\Herm(\A,*)$ with 
$[X,Y]_A$, for any $A \in \Aherm(A,*)$.
\end{description}
\end{lemma}

\Pf
(1) follows from the fact that $(X,Y) \mapsto XAY$ is an associative product, and (2) and (3) from the
fact that $*$ is an antiautomorphism (resp.\ automorphism) of the bracket $[X,Y]_A$.
\EPf
The ``classical Lie algebras'' are obtained by
taking  $\A = M(n,n;\bF)$, the matrix algebra  over  $\bF = \R, \C$ or $\HHH$, with
involution $X^* = \delta(X)^t$ (transposed matrix of $\delta(X)$,  where
$\delta:\bF \to \bF$ is an involution of the base field $\bF$; the various choices for
$\delta$ will be specified below).
For the moment, $(\bF,\delta)$ may be any unital ring with involution. We fix notation and
terminology as follows (cf.\ \cite{BeKi09}):

\msk
\noindent
\begin{tabular}{llll}
family name & label and space & parameter space & Lie bracket 
\cr
\hline
general linear (square)
&$\gL_n(A;\bF):=M(n,n;\bF)$
& $A \in M(n,n;\bF)$
& $[X,Y]_A$
\cr
general linear (rectan.)
&
$\gL_{p,q}(A;\bF):=M(p,q;\bF)$
& $A \in M(q,p;\bF)$
& $[X,Y]_A$
\cr
$(\bF,\delta)$-unitary & $\uu_n(A;\bF,\delta):= \Aherm(n;\bF,\delta)$ & 
$A \in \Herm(n;\bF,\delta)$ & $[X,Y]_A$
\cr
$(\bF,\delta)$-symplectic & $\sp_{n/2}(A;\bF,\delta):= \Herm(n;\bF,\delta)$ & 
$A \in \Aherm(n;\bF,\delta)$ & $[X,Y]_A$
\end{tabular}

\begin{rmk}
The general linear type can be defined, more generally, for {\em associative pairs}
(see \cite{BeKi09}); in the  table above this type appears as {\em rectangular matrices} (second
line). For $p \not=q$ these algebras are never reductive. 
For instance, if $q=p+1$ and $A$ is the matrix having coefficients $1$ on the (first)
diagonal and $0$ else, then $\gL_{p,p+1}(A;\bF)$ is the Lie algebra of the
{\em affine group} of $\bF^p$.
Similarly, for odd $n$, the symplectic type
is never reductive, and we then prefer to call it {\em half-symplectic}.
\end{rmk}

Now we specify the involution of the base field or ring $\bF$.
In the following table, $\bF$ is one of the skew-fields $\R$, $\C$, $\HHH$,
and $\K$ one of the fields $\R$, $\C$. 
Concerning involutions: for
 $\bF=\C$, we always use usual
complex conjugation, and for $\K=\HHH$, if nothing else is specified, we use
``usual'' conjugation $\lambda \mapsto \overline \lambda$
 (minus one in the imaginary part
$\im \HHH$ and one on the center $\R \subset \HHH$).
If we consider $\HHH$ with its  ``split'' involution $\lambda \mapsto
\widetilde \lambda:=j \overline \lambda j^{-1}$, then we write $\widetilde \HHH$. 
For instance,
$\Herm(n;\widetilde \HHH)$ is the space of quaternionic matrices
such that $\widetilde X = X^t$, hence  $\uu_n(1 ;\widetilde \HHH)$
is the Lie algebra in the literature often denoted by $\so^*_{2n}$. 
\msk

\noindent
\begin{tabular}{llll}
family name & label and space & parameter space & Lie bracket 
\cr
\hline
general linear 
&
$\gL_{p,q}(A;\bF):=M(p,q;\bF)$
& $A \in M(q,p;\bF)$
& $[X,Y]_A$
\cr
orthogonal
&
$\oo_n(A;\K):= \Asym(n;\K)$
& $A \in \Sym(n;\K)$
& $[X,Y]_A$
\cr
[half-] symplectic
&$\sP_{n/2}(A;\K):= \Sym(n;\K)$
& $A \in \Asym(n;\K)$
& $[X,Y]_A$
\cr
$\C$-unitary
&$\uu_n(A;\C):= \Aherm(n;\C)$
& $A \in \Herm(n;\C)$
& $[X,Y]_A$
\cr
$\HHH$-unitary
& $\uu_n(A;\HHH):= \Aherm(n;\HHH)$
& $A \in \Herm(n;\HHH)$
& $[X,Y]_A$
\cr
$\HHH$-unitary split
& $\uu_n(A;\widetilde \HHH):= \Aherm(n;\widetilde \HHH)$
& $A \in \Herm(n;\widetilde \HHH)$
& $[X,Y]_A$
\end{tabular}

\msk \nin
Classification up to isomorphy of algebras from the preceding table is easy, by using the
following general 

\begin{lemma} \label{algLemma'}
Let $A,T,S$ be arbitrary elements of the associative algebra $\A$. Then
$$
\g_{TAS} \to \g_A, \quad X \mapsto SXT
$$
is a Lie algebra homomorphism. In particular, if $g,h \in \A^\times$, 
the Lie algebras $\g_{gAh}$ and   $\g_{A}$
are isomorphic, and so are $\uu_{gAg^*}$ and   $\uu_{A}$ (resp.\
$\sp_{gAg^*}$ and   $\sp_{A}$).
\end{lemma}

\Pf 
$[SXT,SYT]_A =SXTASYT-SYTASXT=S[X,Y]_{TAS}T$
\EPf

In the case of matrix algebras over a field, any
$A$ is conjugate to a matrix which is {\em idempotent}, and classification of algebras
is reduced  the well-known classification
of (Hermitian or skew-Hermitian) idempotents. Moreover, if $A$ is idempotent, or more generally
if $A^3 = A$, we may apply  the lemma with $T=S=A$
to get an algebra endomorphism of $\g_A$, which leads to the
{\em fibered structure of homotopes}
to be studied in more detail in Part II.

\section{Classical Lie triple systems}

Next we are going to construct families of ``classical Lie triple systems'' in a similar
way as above. The construction will be based on associative algebras with
{\em several commuting involutions}. However, first of all let us review the preceding
situation (no, or just one, involution) from the point of view of Lie triple systems.

\subsection{Associative algebras with one involution}

 In an associative algebra, we will use the notation
\begin{equation}
T(X,Y,Z):= XYZ + ZYX \, .
\end{equation}

\begin{lemma}\label{LemmaA}
Let $\A$ be an associative algebra, $A \in \A$ and $\alpha(X):=AXA$. Then $\A$ with ternary 
bracket
$$
[X,Y,Z]_A := (XAYAZ + ZAYAX) - (YAXAZ + ZAXAY) = T(X,\alpha Y,Z) - T(Y,\alpha X,Z)
$$
 is the Lie triple system belonging to the Lie algebra
$(\A,[X,Y]_A)$.
\end{lemma}

\Pf
$[[X,Y]_A ,Z]_A = (XAY - YAX)AZ - ZA (XAY - YAX) = T(X,\alpha Y,Z) - T(Y,\alpha X,Z)$.
\EPf

One may note that 
$[X,Y,Z]_{rA} = r^2 [X,Y,Z]_{A}$ for $r \in \K$, whence
$[X,Y,Z]_A = [X,Y,Z]_{-A}$, and, if $\K=\C$, the $c$-dual LTS is obtained in the form
$[X,Y,Z]_{iA} = - [X,Y,Z]_A$.

\begin{lemma}\label{OneInv} \label{LemmaB}
Assume $\tau(X):=X^*$ is an involution of the associative algebra $\A$.
 If $A$ belongs to one of the eigenspaces, then both eigenspaces are
stable under $[X,Y,Z]_A$. This gives four possibilities to combine choices, leading
to the following four families of Lie triple systems:
\begin{description}
\item[(1)] 
Assume $A \in \A^\tau$. Then $\A^{-\tau}$ is the LTS belonging to the unitary Lie algebra $\uu_{A}$, 
and the space $\A^\tau$ with $[X,Y,Z]_A$ is the LTS belongs to the symmetric pair
$
(\g_A , \uu_A) .
$
\item[(2)] 
Assume $A \in \A^{-\tau}$. Then $\A^{\tau}$ is the LTS belonging to the (half-)symplectic
 Lie algebra $\sp_{A}$, 
and the space $\A^{-\tau}$ with $[X,Y,Z]_A$ is the LTS belongs to the symmetric pair
$
(\g_A , \sp_A) .
$
\end{description}
We summarize these statements by the following table:

\ssk \begin{center}
\begin{tabular}{l | cc}
 & $A \in \A^{\tau}$  & $A \in \A^{- \tau }$ 
\cr
\hline
LTS $\A^{\tau}$ & $(\g_{\A,A}, \uu_{\A,A,\tau})$ & $\sp_{A,\tau} $
\cr
LTS $\A^{-\tau}$ & $\uu_{\A,A,\tau}$  & $(\g_{\A,A}, \sp_{\A,A,\tau})$
\cr
\end{tabular}
\end{center}
\end{lemma}

\Pf If $\tau(A)=\mp A$, then $\tau$ is either a Lie algebra automorphism or antiautomorphism of $\g_A$,
and hence in both cases it is a LTS-automorphism, and hence both eigenspaces are
sub-LTS. Now both claims are immediate consequences of the definition of $\sp_A$ and
$\uu_A$. 
\EPf

Note that the LTS of group cases are found on the {\em antidiagonal} of the table; this
reflects the fact that $\tau$ is an {\em anti\/}automorphism of $\A$. 
If one is interested in isomorphism classes of the LTS from the lemma, 
one may observe that the group $\A^\times$ acts on both eigenspaces by
$(a,x) \mapsto ax\tau(a)$, and if $A$ and $A'$ are conjugate under this action,
then the LTS indexed by $A$ and $A'$ are isomorphic.

\subsection{Associative algebras with two commuting involutions}

Now assume that $\A$ carries two commuting involutions $\tau$ and $\tilde \tau$, and let
$\phi := \tau \circ \tilde \tau$. This is an automorphism of order $2$, and conversely,
given an automorphism of order $2$ commuting with $\tau$, we recover $\tilde \tau
= \phi \circ \tau$.
For a given pair $(\tau_1,\tau_2)=(\tau,\tilde \tau)$,
will denote joint eigenspaces by a double superscript, for instance
$$
\A^{(1,-1)}:=
\A^{(\tau, - \tilde \tau)} := \A^\tau \cap \A^{-\tilde \tau} =
\{ X \in \A \mid \, \tau(X)=X, \, \tilde \tau(X) = -X \},
$$
and so on. We thus have a decomposition
$$
\A = \A^{(1,1)}  \oplus \A^{(1, - 1)}  \oplus  \A^{(-1, 1)}  \oplus \A^{(-1, - 1)} ,
$$
and $\A^\phi = \A^{(1,1)}  \oplus  \A^{(-1, - 1)}$ is an associative algebra,
whereas the other spaces are in general not associative algebras. 
Nevertheless, the space $\A^{-\phi}= \A^{(1,-1)}  \oplus  \A^{(-1,1)}$ is an {\em associative
triple system}, i.e., closed under the ternary associative product $XYZ$. In particular, for
$A \in \A^{-\phi}$, the space $\A^{-\phi}$ is stable under the 
associative product $(X,Y) \mapsto XAY$, and hence under the
Lie bracket $[X,Y]_A$. We denote this Lie algebra by $\g_{\A^{-\phi},A}$.

\begin{lemma}\label{TwoInv} 
If $A$ belongs to any one of the joint eigenspaces, 
then all four joint eigenspaces are stable under the triple bracket
$[X,Y,Z]_A$.
\end{lemma}

\Pf
This follows immediately from Lemma \ref{OneInv} applied to $\tau$ and to $\tilde \tau$.
\EPf

Obviously, the lemma generalizes to the case of $k$ commuting involutions: we then
have $2^k$ choices for joint eigenspaces and $2^k$ choices for parameter spaces,
leading to $4^{k}$ families of Lie triple systems. The following theorem shows that in our
case ($k=2$), since the roles of $\tau$ and
$\tilde \tau$ are symmetric, the effective number reduces from $16$ to about $10$:

\begin{thm} \label{TwoInvTh} 
The following holds with respect to the Lie triple bracket $[X,Y,Z]_A$.
\begin{description}
\item[(1)] Assume $A \in \A^{(-\tau,-\tilde \tau)}$. Then the Lie algebras $\sp_{A,\tau}$
and $\sp_{A,\tilde \tau}$ are defined, and their intersection $\h:=\sp_{A,\tau} \cap \sp_{A,\tilde \tau}$
 is  the symplectic Lie algebra
$\sp_{\A^\phi,A,\tau}$ belonging to the involution $\tau$ restricted to the associative algebra
$\A^\phi$. With this notation, we have:
\begin{description}
\item[(i)]
$\A^{(\tau, \tilde \tau)} = \h$ is the LTS belonging to the Lie algebra $\h$; 
\item[(ii)]
$\A^{(-\tau, -\tilde \tau)}$ is the LTS belonging to the symmetric pair
$(\g_{\A^\phi,A}, \h)$;
\item[(iii)]
$\A^{(\tau, -\tilde \tau)}$ is the LTS belonging to the symmetric pair
$(\sp_{A,\tau}, \h)$;
\item[(iv)]
$\A^{(-\tau, \tilde \tau)}$ is the LTS belonging to the symmetric pair
$(\sp_{A,\tilde \tau}, \h)$.
\end{description}
\item[(2)] 
Assume $A \in \A^{(\tau,\tilde \tau)}$. Then the Lie algebras $\uu_{A,\tau}$
and $\uu_{A,\tilde \tau}$ are defined, and their intersection $\h:=\uu_{A,\tau} \cap \uu_{A,\tilde \tau}$
 is  the unitary Lie algebra
$\uu_{\A^\phi,A,\tau}$ belonging to the involution $\tau$ restricted to the associative algebra
$\A^\phi$. With this notation, we have:
\begin{description}
\item[(i)]
$\A^{(-\tau, -\tilde \tau)} = \h$ is the LTS belonging to the Lie algebra $\h$; 
\item[(ii)]
$\A^{(\tau, \tilde \tau)}$ is the LTS belonging to the symmetric pair
$(\g_{\A^\phi,A}, \h)$;
\item[(iii)]
$\A^{(-\tau, \tilde \tau)}$ is the LTS belonging to the symmetric pair
$(\uu_{A,\tau}, \h)$;
\item[(iv)]
$\A^{(\tau, -\tilde \tau)}$ is the LTS belonging to the symmetric pair
$(\uu_{A,\tilde \tau}, \h)$.
\end{description}
\item[(3)] 
Assume $A \in \A^{(\tau,-\tilde \tau)}$. Then $\h:=\uu_{A,\tau} \cap \sp_{A,\tilde \tau}=
\A^{(-\tau, \tilde \tau)}$ is a Lie algebra, and
\begin{description}
\item[(i)]
$\A^{(-\tau, \tilde \tau)}$ is the LTS belonging to the
Lie algebra $\h$;
\item[(ii)]
$\A^{(\tau, -\tilde \tau)}$ is the LTS belonging to the symmetric pair
$(\g_{\A^{-\phi},A}, \h)$;
\item[(iii)]
$\A^{(-\tau, -\tilde \tau)}$ is the LTS belonging to the symmetric pair
$(\uu_{A,\tau}, \h)$;
\item[(iv)]
$\A^{(\tau, \tilde \tau)}$ is the LTS belonging to the symmetric pair
$(\sp_{A,\tilde \tau}, \h)$.
\end{description}
\end{description}
We summarize these statements by the following table:

\ssk \noindent
\begin{tabular}{l | cccc}
 & $A \in \A^{(1,1)}$  & $A \in \A^{(-1,1)}$  & $A \in \A^{(1,-1)}$ & $A \in \A^{(-1,-1)}$ 
\cr
\hline
LTS $\A^{(1,1)}$ & $(\g_{\A^\phi,A}, \uu_{\A^\phi,A,\tau})$ & $(\sp_{A,\tau} , \sp_{A,\tau} \cap \uu_{A,\tilde 
\tau})$
& $( \sp_{A,\tilde \tau} ,\uu_{A,\tau} \cap \sp_{A,\tilde \tau})$ & $\sp_{\A^\phi,A,\tau}$
\cr
LTS $\A^{(-1,1)}$ & $(\uu_{\A,A,\tau}, \uu_{\A^\phi,A,\tau})$ & $(\g_{\A^{-\phi},A}, \sp_{A,\tau} \cap \uu_{A,\tilde 
\tau})$  &$\uu_{A,\tau} \cap \sp_{A,\tilde \tau}$  & $(\sp_{\A,A,\tau}, \sp_{\A^\phi,A,\tau})$
\cr
LTS $\A^{(1,-1)}$ & $(\uu_{\A,A,\tilde \tau}, \uu_{\A^\phi,A,\tau})$ & $\sp_{A,\tau} \cap \uu_{A,\tilde \tau}$ 
& $(\g_{\A^{-\phi},A} ,\uu_{A,\tau} \cap \sp_{A,\tilde \tau})$
& $(\sp_{\A,A,\tilde \tau}, \sp_{\A^\phi,A,\tau})$
\cr
LTS $\A^{(-1,-1)}$ & $\uu_{\A^\phi,A,\tau}$  & $(\uu_{A,\tilde \tau} , \sp_{A,\tau} \cap \uu_{A,\tilde \tau})$ & 
$(\uu_{A,\tau} ,\uu_{A,\tau} \cap \sp_{A,\tilde \tau})$
& $(\g_{\A^\phi,A}, \sp_{\A^\phi,A,\tau})$
\cr
\end{tabular}
\end{thm}

\Pf 
Let us explain the general pattern for $k$ commuting involutions $\tau_1,\ldots,\tau_k$.
Given a vector of eigenvalues $\bs \in \{ \pm 1\}^k$, let $\A^\bs = \cap_{i=1}^k \A^{\bs_i \tau_i}$
be the corresponding joint eigenspace. Fix a vector of eigenvalues $\bt$ and assume that
$A \in \A^\bt$. Then the joint eigenspace $\A^{-\bt}$ is a Lie algebra with respect to the
bracket $[X,Y]_A$ (each $\A^{-t_i \tau_i}$ is a Lie algebra for this bracket, according to 
Lemma \ref{OneInv}, and $\A^{-\bt}$ is the intersection of these algebras), and hence
$\h=\A^{-\bt}$ with $[X,Y,Z]_A$ is the LTS belonging to this Lie algebra. This explains item (i)
in each case (the antidiagonal of the table).  
The other three items correspond to the
symmetric pair given by the direct sum of $\h$ with one of the three joint eigenspaces other than the
one from (i): if $\bs$ is different from $-\bt$, then $\A^\bs$ is not a Lie algebra with respect to
$[X,Y]_A$, but
$\A^\bs \oplus \A^{-\bt}$ is  (this follows since 
$[\A^{s_i \tau_i}, \A^{s_i \tau_i}]_A \subset \A^{-s_i \tau_i}$ if $s_i = t_i$ and 
$[\A^{s_i \tau_i}, \A^{s_i \tau_i}]_A \subset \A^{s_i \tau_i}$ if $s_i = - t_i$).
Moreover, for all $i$  with $s_i = t_i$, the restriction of $\tau_i$ to the Lie algebra
$\A^\bs \oplus \A^{-\bt}$ leads to the same Lie algebra automorphism 
with fixed algebra $\h=\A^{-\bt}$, and hence
the LTS $\A^\bs$ belongs to the symmetric Lie algebra $(\A^\bs \oplus \A^{-\bt},\h)$. 

For $k=2$, the Lie algebras $\A^\bs \oplus \A^{-\bt}$ have explicit descriptions as follows:
let $\bs = \bt$, that is,  $A$ belongs to the underlying space of the LTS in question; then
for $\bt = (-1,-1)$ we get $\A^\bs \oplus \A^{-\bt} = \A^{(-1,-1)} \oplus \A^{(1,1)}$, the fixed
point space of $\phi$ which is an {\em associative} algebra (case (1), (ii));
similarly for $\bt = (1,1)$, whereas for
$\bt = (1,-1)$ we have $\A^\bs \oplus \A^{-\bt} = \A^{(-1,1)} \oplus \A^{(1,-1)}$,
the antifixed space of the associative automorphism $\phi$ (case (3), (ii))
(which corresponds to the associative triple system $\A^{-\phi}$).

If neither $\bs = \bt$ nor $\bs = - \bt$, then  $\A^\bs \oplus \A^{-\bt}$
is equal to one of the spaces $\A^\tau$, $\A^{-\tau}$, $\A^{\tilde \tau}$ or
$\A^{-\tilde \tau}$ with Lie bracket $[X,Y]_A$, leading to the eight remaining cases of
the table. 
\EPf

The presentation in form of a table reveals a remarkable duality between ``space'' (lines)
and ``deformation parameter'' (columns), which is not predicted by the general
theory to be developed in Part II (it reminds Howe's duality of dual pairs
in some respects).
Note that  the diagonal terms in the table are all 
of type ``(general linear, half-symplectic)'' or ``(general linear, unitary)'', whereas the 
antidiagonal terms are algebra cases.

\begin{lemma}\label{strLemma}
If $\tau$ and $\tilde \tau$ are two commuting involutions and $\phi = \tau \circ \tilde \tau$, then 
the group $\Gamma := (\A^\phi)^\times$ acts on all four joint eigenspaces by
$(g,x) \mapsto gx\tau(g)$.  If $A$ and $A'$ are conjugate under this action, then
the corresponding homotope LTS are isomorphic to each other.
\end{lemma}

\Pf
Straightforward calculation (note that $\tau(g)=\tilde \tau(g)$ for $g \in \Gamma$).
\EPf

In all of the following examples, the group $\Gamma$ turns out to be the ``natural''
group acting on the given data, so that the description of its orbits amounts in all
cases to more or less standard results in linear algebra.  This will make
classification up to isomorphy quite easy (in most cases, $A$ will be conjugate
to an idempotent element under this action). However, if we go beyond the standard
examples (for instance, looking at infinite dimensional algebras), then such
a classification is generally completely out of reach. 

\ssk
Theorem \ref{TwoInvTh} contains a great variety of interesting special cases: indeed,
the situation of an associative algebra with two commuting involutions is very common. 
We are going to work out explicitly some of these special cases. 

\begin{rmk}
{\bf (c-duality.)}
In all of the following examples, one may write ``$c$-dual  tables'' in the following way:
consider a real
involutive algebra $(\B,*)$ and complexify it: $\A = \B_\C$, let $\tau$ the $\C$-linear
extension of $*$ and $\tilde \tau(X):=\overline{\tau(X)}$ its $\C$-antilinear extension.
Then the small squares from the preceding table
obtained by taking the middle entries, resp.\ by the ``corner entries'',
will reproduce the tables for $(\B,*)$ from Lemma \ref{OneInv}, whereas
the other eight entries will contain the {\em $c$-dual symmetric pairs} of those.
Similarly, if we complexify a real algebra with two commuting involutions, we get a
complex algebra with three commuting involutions, and the corresponding table
of size $8 \times 8$ will contain together with each LTS also their $c$-duals. 
Replacing $\R$ by $\K$ and $\C$ by $\K[X]/(X^2 +1)$, this construction also applies
over general base rings $\K$. However, for reasons of space we will not write out such
big tables. 
\end{rmk}

\subsection{Case of a matrix algebra with an idempotent}

Let $\A = M(n,n;\K)$, $\tau(X)=X^t$ (transposed matrix) and
$I_{p,q}=\begin{pmatrix} 1_p & 0 \cr 0 & -1_q \end{pmatrix}$.
These data are paradigmatic for the following, slightly more general, situation:
assume given an involutive algebra $(\A,\tau)$ with an idempotent $e$
such that $\tau(e)=e$; then $c:=1 - 2e$ is an element such that $c^2 = 1$, and
$\tilde \tau(x):=c \tau(x) c$ is an involution commuting with $\tau$.

\begin{thm}\label{ProjHomotopeTh} {\rm (Homotopes of projective and of polarized spaces)}
Let $\A = M(n,n;\K)$, fix a decomposition $n=p+q$  and the pair of involutions
$(\tau,\tilde \tau)$ with $\tau(X)=X^t$, $\tilde \tau (X) =I_{p,q}X^t I_{p,q}$, 
so that $\phi(X)=I_{p,q} X I_{p,q}$ is conjugation by $I_{p,q}$.
Then the eigenspaces are
$$
\A^{(1,1)}=\Big\{ \begin{pmatrix} B & 0 \cr 0 & C \end{pmatrix} \mid B \in \Sym(p,\K),C \in \Sym(q,\K) 
\Big\} \cong
\Sym(p,\K) \oplus \Sym(q,\K)
$$
$$
\A^{(1,-1)} = \Big\{ \begin{pmatrix} 0 & A \cr A^t & 0 \end{pmatrix} \mid A \in M(p,q;\K)
\Big\} \cong
M(p,q;\K),
$$
$$
 \A^{(-1,1)} \cong M(q,p;\K), \quad
\A^{(-1,-1)} \cong \Asym(n,\K) \oplus \Asym(n,\K)
$$
and $\A^\phi \cong M(p;\K) \oplus M(q;\K)$, $\A^{-\phi} = M(p,q;\K) \oplus M(q,p;\K)$.
The Lie triple systems from Theorem \ref{TwoInvTh}  are 
explicitly given by the following table (where, for the purpose  of space economy, we write 
symmetric pairs in the form of a quotient):

 \msk \nin
\begin{tabular}{l | cccc}
 & $A=\begin{pmatrix} B & 0 \cr 0 & C \end{pmatrix} \in \A^{(1,1)}$  & $A\in \A^{(-1,1)}$  & $A \in \A^{(1,-1)}$ & $A=\begin{pmatrix} B & 0 \cr 0 & C \end{pmatrix}  \in \A^{(-1,-1)}$ 
\cr
\hline
$\A^{(1,1)} $ & $\frac{ \gL_n(B,\K) \times \gL_n(B,\K)}{\oo_n(B,\K) \times \oo_n(C,\K) }$ 
& $\frac{\sp_{\frac{n}{2}}( \begin{pmatrix} 0 & A \cr -A^t & 0 \end{pmatrix},\K)}{\gL_{p,q}(A,\K)}$
& $\frac{\sp_{\frac{n}{2}}( \begin{pmatrix} 0 & A \cr -A^t & 0 \end{pmatrix},\K)}{\gL_{p,q}(A,\K)}$
& $\sp_{\frac{n}{2}}(B;\K) \times \sp_{\frac{n}{2}}(C;\K)$
\cr
$\A^{(-1,1)}$ & $\frac{ \oo_{n}( \begin{pmatrix} B & 0 \cr 0 & C \end{pmatrix} ,\K)}
{\oo_n(B;\K) \times \oo_n(C;\K)}$  & 
$\gL_{p,q}(A,\K)$  
& $\gL_{p,q}(A,\K)$  
 & $\frac{ \sp_{\frac{n}{2}}( \begin{pmatrix} B & 0 \cr 0 & C \end{pmatrix} ,\K)}
{\sp_{\frac{n}{2}}(B;\K) \times \sp_{\frac{n}{2}}(C;\K)}$ 
\cr
$\A^{(1,-1)}$ &  
$\frac{ \oo_{n}( \begin{pmatrix} B & 0 \cr 0 & -C \end{pmatrix} ,\K)}
{\oo_n(B;\K) \times \oo_n(C;\K)}$  & 
$\gL_{p,q}(A,\K)$  
& $\gL_{p,q}(A,\K)$  
 & $\frac{ \sp_{\frac{n}{2}}( \begin{pmatrix} B & 0 \cr 0 & -C \end{pmatrix} ,\K)}
{\sp_{\frac{n}{2}}(B;\K) \times \sp_{\frac{n}{2}}(C;\K)}$ 
\cr
$\A^{(-1,-1)}$ & $ \oo_p(B,\K) \times \oo_q(C,\K)$  & 
$\frac{\oo_{n}( \begin{pmatrix} 0 & A \cr A^t & 0 \end{pmatrix},\R)}{\gL_{p,q}(A,\K)}$ &  
$\frac{\oo_{n}( \begin{pmatrix} 0 & A \cr A^t & 0 \end{pmatrix},\R)}{\gL_{p,q}(A,\K)}$ &  
 $\frac{ \gL_n(B,\K) \times \gL_n(B,\K)}{\sp_{\frac{p}{2}}(B,\K) \times \sp_{\frac{q}{2}}(C,\K)  }$ 
\end{tabular}
\end{thm}

\Pf
The determination of the eigenspaces is given by straightforward calculations.
Most descriptions of the Lie algebras appearing in the table are also fairly straightforward
from definitions, except perhaps the special form of the ``middle square'':
here, the special feature is that there are two diagonal terms $\gL_{p,q}(A,\K)$ of 
algebra type. This is due to the fact that the Lie algebra $\g_{\A^{-\phi},A}$ is, in our case,
a direct product $\gL_{p,q}(A,\K) \times \gL_{p,q}(A,\K)$ (indeed, for $A$ as in the table,
$\A^{-\phi}$ with $(X,Z)\mapsto XAZ$ is a direct product of associative algebras), 
and hence we get a LTS of
group type. (All this holds, more generally, for associative algebras with idempotent $e$,
as mentioned above, and using the Peirce-decomposition with respect to $e$.)
\EPf

\nin
{\bf Comments.}
For $\K=\R$, on the level of symmetric spaces,
the second and third line of the table describe homotopes of the Grassmannians
$\Gras_p(\R^n)=
\OO(p+q)/\OO(p)\times \OO(q)$ (which arise for $B,C$ being identity matrices).
The first and last  line describe homotopes of certain
 ``polarized symmetric spaces'' (see Section \ref{Sec:Polarized} below):
  notice first that $\A^{(1,1)}$ and $\A^{(-1,-1)}$ have, as
 vector spaces, a natural direct product structure; the corresponding symmetric spaces 
 inherit this product structure, but whereas for the first and last spaces (in the first and last line)
  this product
 structure is {\em global}, for the middle two
 spaces it is only {\em local}: globally,  they are not direct products; for $p=q$
they are homotopes of 
  $\Sp(p,\R)/\Gl(p,\R)$ resp.\ of
$\OO(p,p)/\Gl(p,\R)$ (which are instances of  ``Cayley type symmetric spaces'');
for $p \not= q$ these families never contain reductive
symmetric spaces (and so far seem not to have appeared in the literature).

\subsection{Case of the algebra $\A = M(2,2;\B)$ for an involutive algebra $\B$}

An another important case of application of Theorem \ref{TwoInvTh} is
the algebra of $2 \times 2$-matrices
$\A:=M(2,2;\B)$ with coefficients in $\B$, where
 $\B$ is an associative algebra with involution $*$.
 Then $\A$ is again involutive: let us call
$\tau_1(X):=(X^*)^t$ (``transposed conjugate matrix'') the {\em standard involution of $\A$}.
Besides $\tau_1$, there are at least three other fairly canonical involutions on $\A$.
They are defined using the matrices
\begin{equation}
J:= \begin{pmatrix} 0 & 1 \cr - 1 & 0 \cr \end{pmatrix}, \quad
F:= \begin{pmatrix} 0 & 1 \cr 1 & 0 \cr \end{pmatrix}, \quad
I:= JF = \begin{pmatrix} 1 & 0 \cr 0 & -1 \cr \end{pmatrix} \, .
\end{equation}
For an invertible matrix $B$ let $B_*(X) = BXB^{-1}$ be conjugation by $B$; if
$\tau(B)=\pm B^{-1}$, then $\tau \circ B_* = B_* \circ \tau$ is again an involution. 
Hence we have the following involutions
$$
\tau_J := J_* \circ \tau_1 , \quad
\tau_F:=F_* \circ \tau_1, \quad
\tau_I:= I_* \circ \tau_1 \, .
$$
These involutions commute among each other, and
$\tau_1 \circ \tau_J = J_*$, $\tau_I \circ \tau_F = J_*$, etc.
We have the following explicit formulae:
$$
\tau_1 \begin{pmatrix} a & b \cr c & d \end{pmatrix} =
\begin{pmatrix} a^* & c^* \cr b^* & d^* \end{pmatrix} , \quad
\tau_J \begin{pmatrix} a & b \cr c & d \end{pmatrix} =
\begin{pmatrix} d^* & -b^* \cr -c^* & a^* \end{pmatrix} , 
$$
$$
\tau_F \begin{pmatrix} a & b \cr c & d \end{pmatrix} =
\begin{pmatrix} d^* & b^* \cr c^* & a^* \end{pmatrix} , \quad
\tau_I \begin{pmatrix} a & b \cr c & d \end{pmatrix} =
\begin{pmatrix} a^* & -c^* \cr -b^* & d^* \end{pmatrix} \, .
$$
We call $\tau_J$ the {\em symplectic involution} and $\tau_F$ the
{\em artinian involution of $\A$}.
In the following, we will mainly be interested in the case $\B = M(n,n;\K)$
with $X^* = X^t$ (transposed matrix); then 
$\A = M(2,2;\B) = M(2n,2n;\K)$, and the standard involution $\tau_1$ on
$\A$ is precisely the usual transposed of $2n \times 2n$-matrices. Therefore
the eigenspaces of $\tau_1$ are $\Sym(2n,\K)$ and $\Asym(2n,\K)$.
For the eigenspaces of $\tau_F$, note that $X^t = X$ is equivalent to
$F (FX)^t F = FX$, and hence 
\begin{equation}
\A^{\tau_F} = F \A^{\tau_1} = F \Sym(2n,\K), \quad \A^{-\tau_F} = F \A^{-\tau_1} = F \Asym(2n,\K)
\end{equation}
Similarly, $\A^{\tau_I} = I  \Sym(2n,\K)$ and $\A^{-\tau_I} = I \Asym(2n,\K)$, but
\begin{equation}
\A^{\tau_J} = J \A^{-\tau_1} = J \Asym(2n,\K), \quad \A^{-\tau_J} = J \A^{\tau_1} = J \Sym(2n,\K)
\end{equation}
In the following, we assume that  $\K=\R$.  
Then  $\tau_I \circ \tau_F = J_*$ is conjugation by the standard complex
structure on $\R^{2n}$, and hence its fixed point algebra is $M(n,n;\C)$.

\begin{thm} {\rm (Homotopes of the Siegel half plane)}\label{SiegelTh}
With notation as above, let $\B = M(n,n;\R)$, whence $\A = M(2n,2n;\R)$, and fix the pair of involutions
$(\tau,\tilde \tau):=(\tau_I,\tau_F)$, so that $\phi = J_*$ is conjugation by $J$.
Then the eigenspaces are
$$
\A^{(1,1)} = \Sym(n,\C), \quad \A^{(-1,1)}=I \Herm(n,\C), \quad 
\A^{(1,-1)}=F \Herm(n,\C), \quad \A^{(-1,-1)}=\Asym(n,\C)
$$
and $\A^\phi = M(n,n;\C)$,
and the Lie triple systems from Theorem \ref{TwoInvTh}  are 
explicitly given by the following table:

 \msk \nin
\begin{tabular}{l | cccc}
 & $A \in \A^{(1,1)}$  & $A \in \A^{(-1,1)}$  & $A \in \A^{(1,-1)}$ & $A \in \A^{(-1,-1)}$ 
\cr
\hline
$\A^{(1,1)} $ & $(\gL_n(A,\C) ,  \oo_n(A,\C))$ & $(\sp_n(IA,\R) , \uu_n(IA,\C))$
& $(\sp_n(FA,\R) , \uu_n(FA,\C))$ & $\sp_{\frac{n}{2}}(A,\C)$
\cr
$\A^{(-1,1)}$ & $(\oo_{2n}(IA,\R), \oo_n(A,\C)$ & 
$(\gL_n(IA,\C),\uu_n(IA,\C))$ & $\uu_n(FA,\C)$ & $(\sp_n(FA,\R),\sp_{\frac{n}{2}}(A,\C))$
\cr
$\A^{(1,-1)}$ &  $(\oo_{2n}(FA,\R), \oo_n(A,\C)$ &  $\uu_n(IA,\C)$ & 
$(\gL_n(FA,\C),\uu_n(FA,\C))$ & $(\sp_n(IA,\R),\sp_{\frac{n}{2}}(A,\C))$
\cr
$\A^{(-1,-1)}$ & $\oo_n(A,\C)$ & $(\oo_{2n}(FA,\R),\uu_n(IA,\C))$ &  
$(\oo_{2n}(IA,\R),\uu_n(FA,\C))$ & $(\gL_n(A,\C),\sp_{\frac{n}{2}}(A,\C))$
\end{tabular}
\end{thm}

\Pf
First of all, we describe the eigenspaces: as noticed above, $\A^\phi =M(n,n;\C)$, that is, a complex matrix $a + i b$ will be identified with
$X = \begin{pmatrix} a & b \cr - b & a \end{pmatrix} \in \A$.
Similarly,
\begin{equation}
\A^{-\phi} = \Big\{  \begin{pmatrix} a & b \cr  b & - a \end{pmatrix} \mid \, a, b \in M(n,n;\R) \Big\} =
I \,  M(n,n;\C)
\end{equation}
is the space of $\C$-antilinear operators (whose base point $I$ is complex conjugation). 
As an associative triple system, it is isomorphic to $M(n,n;\C)$. 
Now note that
the complex matrix $a + ib$ is symmetric iff $a$ and $b$ are real symmetric matrices, which
means that $X$ is fixed under the involution $\tau_I$ (and not under the standard involution!). 
Summing up,
\begin{equation}
\Sym(n,\C) = \A^{J_*} \cap \A^{\tau_I} =
\A^{\tau_F} \cap \A^{\tau_I} = \A^{(1,1)} \, .
\end{equation}
Similarly 
$$
\A^{(-1,-1)} = \A^{J_*} \cap \A^{-\tau_I} =
\big\{ \begin{pmatrix} a & b \cr -b & a \end{pmatrix} \mid \, a^t = -a, b^t = -b \big\}= \Asym(n,\C).
$$
Next,
 $a+ib$ is Hermitian iff $a=a^t$ and $b=-b^t$, whence $\tau_1(A)=A$ and thus
\begin{eqnarray}
\A^{\tau_J} \cap \A^{\tau_1} & =  & 
 \A^{J_*} \cap \A^{\tau_1}=M(n,n;\C) \cap \Sym(2n;\R) = \Herm(n,\C)
\cr
\A^{(\tau,-\tilde \tau)} &=& \A^{-J_*} \cap \A^{\tau_I} = 
\big\{ \begin{pmatrix} a & b \cr b & -a \end{pmatrix} \mid \, a^t = -a, b^t = b \big\} =
I \, \Herm(n,\C)
\cr
\A^{(-\tau,\tilde \tau)} &=& \A^{-J_*} \cap \A^{-\tau_I} = 
\big\{ \begin{pmatrix} a & b \cr b & -a \end{pmatrix} \mid \, a^t = a, b^t = -b \big\} =
I  \, \Aherm(n,\C) = F  \, \Herm(n,\C) \, .
\end{eqnarray}

Now, for $A$ belonging to a joint eigenspace, let us determine the Lie triple systems.
The Lie triple structure for the four ``corners'' of the table is simply the one of Lemma
\ref{OneInv}, applied to the algebra $\A^\phi = M(n,n;\C)$ with involution being usual transposed.
Similarly, the four ``inner entries'' of the table correspond to the associative pair
$\A^{-\phi} = I M(n,n;\C)$ with involution corresponding to $X \mapsto \overline X^t$ on $M(n,n;\C)$.
For the remaining eight Lie triple systems $(\g,\h)$, the stabilizer algebra $\h$ is given by
the algebra on the intersection of the same column with the diagonal. 
The Lie algebra $\g$ is then one of the algebras
$\A^\tau = I \Sym(2n,\R)$, $\A^{\tilde \tau} = F \Sym(2n,\R)$ (symplectic) or
$\A^{-\tau} = I \Asym(2n,\R)$, $\A^{-\tilde \tau} = F \Asym(2n,\R)$ (orthogonal).
One just has to pay attention that all objects are well-defined; for instance,
if $A \in I \Herm(n,\C)$, 
then $IA \in \Herm(n,\C)$ is a real symmetric matrix (hence $\oo_{2n}(IA;\R)$ is well-defined), and 
$JI A \in J \Herm(n,\C) = \Aherm(n,\C)$ is skew-Hermitian (hence $\uu_n(FA;\C)$ is well-defined)
and also skew-symmetric (hence $\sp_n(FA,\R)$ is well-defined),
and so on.
\EPf

\nin {\bf Comments.} 
In the preceding theorem,
 $\A^{(1,1)}$ and $\A^{(-1,-1)}$ have a natural underlying structure of {\em complex}
vector space,  hence on the level of symmetric spaces we obtain homotopes of {\em (pseudo-)
Hermitian symmetric spaces}. In particular,
the first line describes homotopes of the Siegel half plane $\Sp_n(\R)/\UU(n)$. Indeed,
the third family of symmetric pairs in this line can be written in the form
$(\sp_{n}  \big( \begin{pmatrix} a & b \cr -b & a \end{pmatrix} , \R \big) ,\uu_n(a+ib,\C))$ with
$a \in \Asym(n,\R)$ and $b \in \Sym(n,\R)$, so that the Siegel half plane corresponds to 
the choice $a=0$, $b=1$. 

\subsection{Algebras of quaternionic matrices}

The next two theorems deal with quaternionic matrices. There are two different viewpoints:
we may consider $M(n,n;\HHH)$ as a real form of the complex algebra $\A = M(2n,2n;\C)$,
thus thinking of quaternionic matrices as a special kind of complex matrices; or we may
work intrinsically with matrices having coefficients in $\HHH$.
We start with the latter viewpoint: recall the two involutions of $\HHH$ from Chapter 2, and
let $\A = M(n,n;\HHH)$ and consider the two involutions
$\tau(X):=\overline X^t$ and $\tilde \tau(X)=\tilde X^t$. 
They commute, and $\phi=\tau \circ \tilde \tau$ is the automorphism acting on each coefficient
by conjugation with the quationion $j$.
The field fixed under conjugation by $j$ in $\HHH$ is $\R \oplus j \R \cong \C$, and
hence $\A^\phi \cong M(n,n;\C)$.
Note that $\widetilde{(j a_{ij})} =- \widetilde a_{ij} j = -j \overline a_{ij}$, which implies
\begin{equation}\label{HermQuat}
j \Herm(n,\HHH) = \Aherm(n,\widetilde \HHH), \quad
j \Aherm(n,\HHH)= \Herm(n,\widetilde \HHH) \, .
\end{equation}

\begin{thm} {\rm (Homotopes of quaternionic type. I)}\label{QuaternionicTh1}
With notation as above, let $\A = M(n,n;\HHH)$ and fix the pair of involutions
$(\tau,\tilde \tau)$.
Then the eigenspaces are
$$
\A^{(1,1)} = \Herm(n,\C), \quad \A^{(-1,1)} = i \Sym(n,\C), \quad 
\A^{(1,-1)} =i \Asym(n,\C), \quad \A^{(-1,-1)}=\Aherm(n,\C)
$$
and $\A^\phi \cong M(n,n;\C)$,
and the Lie triple systems from Theorem \ref{TwoInvTh}  are 
explicitly given by the following table

 \msk \nin
\begin{tabular}{l | cccc}
 & $A \in \A^{(1,1)}$  & $A \in \A^{(-1,1)}$  & $A \in \A^{(1,-1)}$ & $A \in \A^{(-1,-1)}$ 
\cr
\hline
$\A^{(1,1)} $ & $(\gL_n(A,\C) ,  \uu_n(A,\C))$ & $(\uu_n(A,\widetilde \HHH) , \sp_{\frac{n}{2}}(A,\C))$
& $(\uu_n(A,\HHH) , \oo_n(A,\C))$ & $\uu_n(A,\C)$
\cr
$\A^{(-1,1)}$ & $(\uu_{n}(A,\widetilde \HHH), \uu_n(A,\C))$ & 
$(\gL_n(A,\C),\sp_{\frac{n}{2}}(A,\C))$ & $\oo_n(A,\C)$ & $(\uu_n(A,\widetilde \HHH),\uu_n(A,\C))$
\cr
$\A^{(1,-1)}$ &  $(\uu_{n}(A,\HHH), \uu_n(A,\C))$ &  $\sp_{\frac{n}{2}}(A,\C)$ & 
$(\gL_n(A,\C),\sp_{\frac{n}{2}}(A,\C))$ & $(\uu_n(A,\HHH),\uu_n(A,\C))$
\cr
$\A^{(-1,-1)}$ & $\uu_n(A,\C)$ & $(\uu_{n}(A,\widetilde \HHH),\sp_{\frac{n}{2}}(A,\C))$ &  
$(\uu_{n}(A,\HHH),\oo_n(A,\C))$ & $(\gL_n(A,\C),\uu_n(A,\C))$
\end{tabular}
\end{thm}

\Pf 
For the eigenspaces,
let $A = (a_{ij}) \in M(n,n;\HHH)$. Then $A \in \A^{(1,1)}$ iff $A \in \A^\phi$
and $\overline A^t = A$, iff $a_{ij} \in \R \oplus j \R = \C$ and
$a_{ij} = \overline a_{ji}$ in $\C$, that is, iff $A \in \Herm(n,\C)$. 
Similarly, $A \in \A^{(1,-1)}$ iff $a_{ij} \in i \R \oplus k \R$
and $a_{ij} = - a_{ji}$, that is, iff $iA$ is a complex symmetric matrix. 
The remaining computations are similar as above, so we omit details.
%
\EPf

\nin {\bf Comments.}
Notice that here the second and third line correspond to homotopes of (pseudo-) Hermitian
symmetric spaces. Indeed, by direct inspection we see that
the diagram contains precisely the $c$-dual symmetric pairs from those given in Theorem
\ref{SiegelTh}. In particular, the third line contains homotopes of the Siegel half plane 
(or, equivalently, of its compact dual $\Sp(n)/\UU(n)$).

\begin{thm} {\rm (Homotopes of quaternionic type. II)}\label{QuaternionicTh2}
Let $\A = M(2n,2n;\C)$ and fix the pair of involutions
$(\tau,\tilde \tau)$ with $\tau(X)=I X^t I$ and $\tilde \tau(X)=F \overline X^t F$.
Then the eigenspaces are
$$
\A^{(1,1)} \cong \Aherm(n,\HHH), \quad \A^{(-1,1)} \cong \Herm(n,\HHH), \quad 
\A^{(1,-1)} \cong \Aherm(n,\HHH), \quad \A^{(-1,-1)}=\Herm(n,\HHH)
$$
and $\A^\phi \cong M(n,n;\HHH)$,
and the Lie triple systems from Theorem \ref{TwoInvTh}  are 
explicitly given by the following table

 \msk \nin
\begin{tabular}{l | cccc}
 & $A \in \A^{(1,1)}$  & $A \in \A^{(-1,1)}$  & $A \in \A^{(1,-1)}$ & $A \in \A^{(-1,-1)}$ 
\cr
\hline
$\A^{(1,1)} $ & $(\gL_n(A,\HHH) ,  \uu_n(A,\widetilde \HHH))$ & 
$(\uu_{2n}(A,\C) , \uu_n(A,\HHH ))$
& $(\sp_n(A,\C) , \uu_n(A,\widetilde \HHH))$ & $\uu_n(A, \HHH)$
\cr
$\A^{(-1,1)}$ & $(\uu_{2n}(A,\C), \uu_n(A,\widetilde \HHH))$ & 
$(\gL_n(A,\HHH),\uu_n(A, \HHH))$ & $\uu_n(A,\widetilde \HHH)$ & 
$(\oo_{2n}(A,\C),\uu_n(A, \HHH))$
\cr
$\A^{(1,-1)}$ &  $(\oo_{2n}(A,\C), \uu_n(A,\widetilde \HHH))$ &  $\uu_n(A, \HHH)$ & 
$(\gL_n(A,\HHH),\uu_n(A,\widetilde \HHH))$ & $(\uu_{2n}(A,\C),\uu_n(A, \HHH))$
\cr
$\A^{(-1,-1)}$ &
$\uu_n(A,\widetilde \HHH)$ & $(\sp_{n}(A,\C),\uu_n(A,\HHH))$ &  
$(\uu_{2n}(A,\C),\uu_n(A,\widetilde \HHH))$ & $(\gL_n(A,\HHH),\uu_n(A, \HHH))$
\end{tabular}
\end{thm}

\Pf
Since $\phi(X)=IF \overline X FI = J \overline X J^{-1}$, the algebra $\A^\phi$ is the
real form $M(n,n;\HHH)$ of $M(2n,2n;\C)$, and $\A^{-\phi}=i M(n,n;\H)$ is, as associative
triple system, again isomorphic to $M(n,n;\HHH)$.
The restriction of $\tau$ to $\A^\phi = M(n,n;\HHH)$ has the same effect as the involution
$X \mapsto \tilde X^t$ there, thus $\A^{(1,1)} =  \Herm(n,\tilde \HHH) \cong \Aherm(n,\HHH)$
(mind the isomorphism (\ref{HermQuat})), and similarly for the other joint eigenspaces.
This gives the inner  $2 \times 2$-square and the square formed by the corner entries.
For the remaining eight spaces, observe that
$\A^{-\tau} \cong \oo_{2n}(A,\C)$ and
$\A^{-\tilde \tau} \cong \uu_{2n}(A,\C)$; combining with the known diagonal entries, this
permets to complete the table.
\EPf

\section{Homotopes of classical real symmetric spaces}

\subsection{Classical Groups}

The following lemma leads to a definition
of Lie groups and algebraic groups
corresponding to the algebras defined in Section \ref{ClassLieAlgebras}:

\begin{lemma} \label{groupLemma}
Let $\A$ be an associative algebra with involution $*$ and fix $A \in \A$.
\begin{description}
\item[(1)]  
The  product
 $X \cdot_A Y = X + Y - XAY$ defines a group structure on the set
$$
G_A:=\{ X \in \A \mid 1 - XA \in \A^\times \} .
$$ 
The neutral element is $0$, and
the inverse is  $j_A(X)= - (1 - XA)^{-1} X$.
More generally, for any associative pair $(\A^+,\A^-)$ the same formulae define a group
structure on $\A^+$, for each $A \in \A^-$.
\item[(2)]
If $A^* = A$, then the following is a subgroup of $G_A$:
$$
U_A:= \{ X \in G_A \mid j_A(X) = X^* \} =
\{ X \in G_A \mid \, X^* + X = X^* A X \}
$$
\item[(3)]
If $A^* = -A$, then the following is a subgroup of $G_A$:
$$
S_A:= \{ X \in G_A \mid j_A(X) = - X^* \} =
\{ X \in G_A \mid \, X^* - X = X^* A X \}
$$
\item[(4)]
If $\A$ is finite dimensional over $\K=\R$, then these groups are Lie groups having as
Lie algebra the corresponding algebra from Lemma \ref{algLemma}.
\end{description}
\end{lemma}

\Pf
(1) Associativity is checked by direct computation; for inversion check first that
$X \mapsto 1-AX$ is a homomorphism from $\cdot_A$ to the usual product. 
For the  case of an associative pair, see \cite{BeKi09}. 

(2) If $A^*=A$, then $*$ is an antiautomorphism of $G_A$. Note  that the
condition $j_A(X)=X^*$ is equivalent to $-X = (1-XA)X^*$, hence to
$X^* + X = XAX^*$. 
This proves (2), and (3) is shown similarly.

(4) This follows easily by differentiating  (see \cite{BeKi09} for a more algebraic argument).
\EPf

\nin
If $\bF$ is a base ring with involution $\delta$ and $\A = M(n,n;\bF)$ we write also 

\msk
\noindent
\begin{tabular}{llll}
label &  underlying set  & parameter space & product \cr
\hline
 $\UU_n(A;\bF,\delta)$ & $:=\{ X \in \Gl_n(A;\bF) |  X + \delta(X)^t = \delta(X)^t AX \}$
& $A \in \Herm(n;\bF,\delta)$
& $X \cdot_A Y$
\cr
 $\Sp_n(A;\bF,\delta)$ & $:=\{ X \in \Gl_n(A;\bF) |  X - \delta(X)^t = \delta(X)^t AX \}$
& $A \in \Aherm(n;\bF,\delta)$
& $X \cdot_A Y$
\end{tabular}

\msk \nin
and using the standard involutions of $\R,\C,\HHH$  leads to the following table of classical groups

\msk
\noindent
\begin{tabular}{llll}
label &  underlying set  & parameter space & product \cr
\hline
 $\Gl_n(A;\bF)$ & $:=\{ X \in M(n,n;\bF) | 1 - AX \, \mbox{invertible} \}$
& $A \in M(n,n;\bF)$
& $X \cdot_A Y$
\cr
$\Gl_{p,q}(A;\bF)$ & $:=\{ X \in M(n,n;\bF) | 1 - AX \, \mbox{invertible} \}$
& $A \in M(q,p;\bF)$
& $X \cdot_A Y$
\cr
$\OO_n(A;\K)$ & $:= \{ X \in \Gl_n(A,\K) |  X + X^t = X^t AX \}$
& $A \in \Sym(n;\K)$
& $X \cdot_A Y$
\cr
$\Sp_{n/2}(A;\K)$ & $:= \{ X \in \Gl_n(A,\K) |  X - X^t = X^t AX \}$
& $A \in \Asym(n;\K)$
& $X \cdot_A Y$
\cr
 $\UU_n(A;\C)$ & 
$:= \{ X \in \Gl_n(A,\K) |  X + \overline X^t = \overline X^t AX 
\}$
& $A \in \Herm(n;\C)$
& $X \cdot_A Y$
\cr
 $\UU_n(A;\HHH)$ & $:= \{ X \in \Gl_n(A,\HHH) |  
X + \overline X^t = \overline X^t AX 
\}$
& $A \in \Herm(n;\HHH)$
& $X \cdot_A Y$
\cr
 $\UU_n(A;\widetilde \HHH)$ & $:= \{ X \in \Gl_n(A,\HHH) |  
X + \widetilde X^t = \widetilde X^t AX 
\}$
& $A \in \Herm(n;\widetilde \HHH)$
& $X \cdot_A Y$
\end{tabular}

\msk \noindent
Concerning classification up to isomorphy, the same remarks as in Section 2 hold:
for all $g,h \in \Gl(n,\K)$, the groups
$\Gl_n(gAh;\K)$ and $\Gl_n(A;\K)$ are isomorphic, via the map
$X \mapsto gXh$. In particular,  $\Gl_n(A;\K)$ and
 $\Gl_n(-A;\K)$ are isomorphic. 
 Similarly,
 unitary or (half-)symplectic groups labelled with
 $A$ and $A' = g A \tau(g)$ for $g \in  \Gl(n,\bF)$ are isomorphic.

\subsection{Classical Symmetric Spaces}

Now we are ready to give a list of classical symmetric spaces $G/H$ and their homotopes.
The labelling given below corresponds to the classification of Jordan triple systems
(see Part II; at this point the labelling may look rather inconsequent -- in fact, it
follows the one from \cite{Be00}, Chapter XII;  in particular,
1.1, 1.2, 1.3 are real forms of the complex type 1, and so on). The term ``classical symmetric space'' is used here
for symmetric spaces corresponding to the matrix families of Jordan triple systems.
Exceptional spaces and the ``semi-exceptional'' family of spin factors are not considered here.

\begin{thm}\label{ClassificationTheorem}
The following tables contain symmetric spaces $M_\alpha = G_\alpha/H_\alpha$ that
are homotopes of each other in the following sense:
fix an underlying real vector space $V^+$; in case
$V^+$ is a space $M(p,q;\bF)$ of rectangular matrices, we let $V^- := M(q,p;\bF)$,  and 
 in all other cases (spaces of symmetric, Hermitian or skew-Hermitian
matrices) we let $V^- := V^+$. 
For each such pair $(V^+,V^-)$ of vector spaces,
define families of linear maps $\alpha:V^+ \to V^-$ as in the tables;
then $V^+$ with triple bracket
$$
[X,Y,Z]_\alpha := T(X,\alpha Y,Z) - T(Y,\alpha X,Z),
$$
where $T(X,Y,Z)=XYZ + ZYX$ (with $XYZ$ and $ZYX$ being usual matrix products),
is a Lie triple system, and it is the LTS belonging to the symmetric space
$G_\alpha/H_\alpha$ in the corresponding line of the following tables: 
\end{thm}

\nin
{\bf
Spaces of rectangular matrices}

\msk
\nin
{\bf 1.} $V^+ = M(p,q;\K)$, $V^-=M(q,p;\K)$, $\K=\R,\C$:

\msk
\noindent
\begin{tabular}{llll}
label & symmetric space $G_\alpha/H_\alpha$ & $\alpha:V^+ \to V^- $ &  parameter set 
 \cr
\hline
1.a & group case $\Gl_{p,q}(A,\K)$  & $\alpha(X)=AXA$ &  $A \in M(q,p;\K)$
\cr
1.a' & $\Gl_{p,q}(A,\K[i])/\Gl_{p,q}(A,\K)$ & $\alpha(X)=-AXA$ &  $A \in M(q,p;\K)$
\cr
1.b & $\OO_{p+q}( \begin{pmatrix} A & 0 \cr 0 & B \end{pmatrix};\K) /
\OO_p(A;\K) \times \OO_q(B;\K)$ & $\alpha(X)=AX^t B$ 
 & $A \in \Sym(p,\K), B \in \Sym(q,\K)$
\cr
1.c & $\Sp_{\frac{p+q}{2}}( \begin{pmatrix} A & 0 \cr 0 & B \end{pmatrix};\K) /
\Sp_{\frac{p}{2}}(A;\K) \times \Sp_{\frac{q}{2}}(B;\K)$  & $\alpha(X)=AX^t B$ 
& $A \in \Asym(p,\K), B \in \Asym(q,\K)$
\end{tabular}

\msk \nin
{\bf 1. cases of $\C$-antilinear $\alpha$:} $V^+ = M(p,q;\C)$, $V^-=M(q,p;\C)$: 

\msk
\noindent
\begin{tabular}{llll}
label & symmetric space $G_\alpha/H_\alpha$ & $\alpha:V^+ \to V^-$ & parameter set 
 \cr
\hline
1.A & $\Gl_{p,q}(A;M(2,2;\R))/ \Gl_{p,q}(A;\C)$ & $\alpha(X)=A \overline{XA}$ & $A \in M(q,p;\C)$
\cr 
1.A' & $\Gl_{p,q}(A;\HHH)/ \Gl_{p,q}(A;\C)$  & $\alpha(X)=-A \overline{XA}$ & $A \in M(q,p;\C)$
\cr
1.B &  $\UU_{p+q}( \begin{pmatrix} A & 0 \cr 0 & B \end{pmatrix};\C) /
\UU_p(A;\C) \times \UU_q(B;\C)$ & $\alpha(X)=A\overline X^t B$ &
 $A \in \Herm(p,\C), B \in \Herm(q,\C)$
\end{tabular}

\msk
\noindent
{\bf 1.3} $V^+ = M(p,q;\HHH)$, $V^- =M(q,p;\HHH)$

\msk
\noindent
\begin{tabular}{llll}
label & symmetric space $G_\alpha/H_\alpha$ & $\alpha:V^+ \to V^-$ &  parameter set 
 \cr
\hline
1.3.a & group case $\Gl_{p,q}(A,\HHH)$   & $\alpha(X)=AXA$ &  $A \in M(q,p;\HHH)$
\cr
1.3.a' &  $\Gl_{p,q}(A,M(2,2;\C)) / \Gl_{p,q}(A,\HHH)$ & $\alpha(X)=-AXA$ & $A \in M(q,p;\HHH)$
\cr
1.3.b & $\UU_{p+q}( \begin{pmatrix} A & 0 \cr 0 & B \end{pmatrix};\HHH) /
\UU_p(A;\HHH) \times \UU_q(B;\HHH)$ & 
$\alpha(X)=A\overline X^t B$ &
$A \in \Herm(p,\HHH),  B \in \Herm(q,\HHH)$
\cr
1.3.c & $\UU_{p+q}( \begin{pmatrix} A & 0 \cr 0 & B \end{pmatrix};\widetilde \HHH) /
\UU_p(A;\widetilde \HHH) \times \UU_q(B;\widetilde  \HHH)$
 & $\alpha(X)=A \widetilde X^t B$ &
 $A \in \Herm(p,\widetilde \HHH),  B \in \Herm(q,\widetilde  \HHH)$
\end{tabular}

\msk \nin
{\bf Spaces of symmetric matrices}

\msk
\noindent
{\bf 2.} $V:=V^+ = V^- = \Sym(n,\K)$, $\K = \R,\C$:

\msk
\noindent
\begin{tabular}{llll}
label & symmetric space $G_\alpha/H_\alpha$ & $\alpha:V \to V$   &  parameter set 
 \cr
\hline
2.a & $\Gl_n(A;\K)/\OO_n(A;\K)$ & $\alpha(X)=AXA$ &  $A \in \Sym(n,\K)$ 
\cr
2.a' & $\UU_n(A;\K[i])/\OO_n(A;\K)$ &  $\alpha(X)=-AXA$ &  $A \in \Sym(n,\K)$ 
\cr
2.b & group space $\Sp_{\frac{n}{2}}(A;\K)$ &  $\alpha(X)=AXA$ &  $A \in \Asym(n,\K)$ 
\cr
2.b' & $\Sp_{\frac{n}{2}}(A;\K[i])/\Sp_{\frac{n}{2}}(A;\K)$  & $\alpha(X)=-AXA$ 
  & $A \in \Asym(n,\K)$
\end{tabular}

\msk
\noindent
{\bf 2. case of antilinear $\alpha$:} $V=V^+ = V^- = \Sym(n,\C)$:

\msk
\noindent
\begin{tabular}{llll}
label & symmetric space  $G_\alpha/H_\alpha$ & $\alpha:V \to V$  &  parameter set 
 \cr
\hline
2.A & $\UU_n(A;\HHH)/\UU_n(A; \C)$   & $\alpha(X)= A \overline X A$ & $A \in \Herm(n,\C)$ 
\cr
2.A' & $\Sp_n \big( \begin{pmatrix} b & a \cr -a & b \end{pmatrix} \big) / \UU_n(b+ia,\C)$ 
& $\alpha(X)= -A \overline X A$ & $A = a + ib \in \Herm(n,\C)$ 
\end{tabular}


\msk \nin
{\bf Spaces of skewsymmetric matrices}

\msk
\noindent
{\bf 3.} $V=V^+ = V^- = \Asym(n,\K)$, $\K = \R,\C$:

\msk
\noindent
\begin{tabular}{llll}
label & symmetric space  $G_\alpha/H_\alpha$ & $\alpha:V \to V$  &  parameter set 
 \cr
\hline
3.a & $\Gl_n(A;\K)/\Sp_{n/2}(A;\K)$ & $\alpha(X)=AXA$ & $A \in \Asym(n,\K)$ 
\cr
3.a' & $\UU_n(A;\K[i])/\Sp_{n/2}(A;\K)$ & $\alpha(X)=-AXA$  & $A \in \Asym(n,\K)$ 
\cr
3.b & group case $\OO_n(A;\K)$ & $\alpha(X)=AXA$  & $A \in \Sym(n,\K)$ 
\cr
3.b' & $\OO_n(A;\K[i])/\OO_n(A;\K)$ & $\alpha(X)=-AXA$  & $A \in \Sym(n,\K)$
\end{tabular}

\msk
\noindent
{\bf 3.  Case of antilinear $\alpha$:} $V=V^+ = V^- = \Asym(n,\C)$

\msk
\noindent
\begin{tabular}{llll}
label & symmetric space $G_\alpha/H_\alpha$ & $\alpha:V \to V$   &  parameter set 
 \cr
\hline
3.A & $\UU_n(A,\widetilde \HHH)/\UU_n(A,\C)$   & $\alpha(X)=A\overline XA$  & $A \in i\Herm(n,\C)$ 
\cr
3.A' & $\OO_{2n} \big( \begin{pmatrix} a & b \cr -b & a \end{pmatrix},\R \big) / \UU_n(b+ia,\C)$ 
& $\alpha(X)=-A \overline XA$ 
 & $A = a + ib \in \Herm(n,\C)$ 
\end{tabular}


\msk \nin
{\bf Spaces of Hermitian matrices}

\msk
\noindent
{\bf 1.1}  $V=V^+ = V^- = \Herm(n,\C)$

\msk
\noindent
\begin{tabular}{llll}
label & symmetric space $G_\alpha/H_\alpha$ & $\alpha:V \to V$  &  parameter set 
 \cr
\hline
1.1.a & $\Gl_n(A,\C)/\UU_n(A,\C)$ & $\alpha(X)=AXA$ &  $A \in \Herm(n,\C)$ 
\cr
1.1.a' & group case $\UU_n(A,\C)$ & $\alpha(X)=-AXA$ & $A \in \Herm(n,\C)$ 
\cr
1.1.b & $\UU_n(A,\widetilde \HHH)/\OO_n(A,\C)$ 
& $\alpha(X)=A \overline{XA}^t$ & $A \in \Sym(n,\C)$ 
\cr
1.1.b' & $\OO_{2n}( \begin{pmatrix} a & b \cr b & -a \end{pmatrix};\R) /
\OO_n(a+ib;\C)$   
& $\alpha(X)=-A \overline{XA}^t$ & $A = a + ib \in \Sym(n,\C)$ 
\cr
1.1.c & $\Sp_{n}(\begin{pmatrix} a & b \cr b & -a \end{pmatrix};\R) /
\Sp_{n/2}(a+ib;\C)$ & $\alpha(X)=A \overline{XA}^t$   & $A = a + ib  \in \Asym(n,\C)$ 
\cr
1.1.c' & $\UU_n(A,\HHH)/\Sp_{n/2}(A,\C)$  
& $\alpha(X)=-A \overline{XA}^t$ & $A \in \Asym(n,\C)$ 
\end{tabular}

\msk
\noindent
{\bf 3.1}  $V^+ = V^- = \Herm(n,\HHH)$:

\msk
\noindent
\begin{tabular}{llll}
label & symmetric space $G_\alpha/H_\alpha$ & $\alpha:V \to V$  &  parameter set 
 \cr
\hline
3.1.a & $\Gl_n(A,\HHH)/\UU_n(A,\HHH)$ & $\alpha(X)=A XA$ & $A \in \Herm(n,\HHH)$ 
\cr
3.1.a' & $\UU_{2n}(IA,\C)/\UU_n(A,\HHH)$  & $\alpha(X)=-A XA$ &$A \in \Herm(n,\HHH)$ 
\cr
3.1.b & group case $\UU_n(A,\widetilde \HHH)$ & $\alpha(X)=A \overline{XA}$ & $A \in \Herm(n,\HHH)$ 
\cr
3.1.b' & $\OO_{2n}(IA,\C)/\UU_n(A,\widetilde \HHH)$ & $\alpha(X)=-A \overline{XA}$
& $A \in \Herm(n,\HHH)$ 
\end{tabular}

\msk
\noindent
{\bf 2.2}  $V^+ = V^- = \Herm(n,\widetilde \HHH)$:

\msk
\noindent
\begin{tabular}{llll}
label & symmetric space $G_\alpha/H_\alpha$ & $\alpha:V \to V$  &  parameter set 
 \cr
\hline
2.2.a & $\Gl_n(A,\HHH)/\UU_n(A,\widetilde \HHH)$ & $\alpha(X)=A XA$ 
& $A \in \Herm(n,\widetilde \HHH)$ 
\cr
2.2.a' & $\UU_{2n}(IA,\C)/\UU_n(A,\widetilde \HHH)$ & $\alpha(X)=-A XA$ 
& $A \in \Herm(n,\tilde \HHH)$ 
\cr
2.2.b & group case $\UU_n(A,\HHH)$ & $\alpha(X)=A \widetilde X \overline A$ 
& $A \in \Herm(n,\widetilde \HHH)$ 
\cr
2.2.b' & $\Sp_{2n}(IA,\C)/\UU_n(A, \HHH)$ & $\alpha(X)=-A \widetilde X \overline A$ 
& $A \in \Herm(n,\widetilde \HHH)$ 
\end{tabular}

\msk \nin
{\em Proof of the theorem.}
All descriptions arise from Theorem \ref{TwoInvTh} by identifying a joint eigenspace 
$\A^{(i,j)}$ with a matrix space $V^+$. Technical complications arise by the fact that
a given matrix space $V^+$ may be realized in several different ways as a joint eigenspace by
using different algebras; this then leads to the various formulas for the endomorphisms
$\alpha$, involving, besides the expression $AXA$ occuring already in Lemma \ref{LemmaA},
transposition and complex or quaternionic conjugation of matrices.
When carrying out the computations, one may mind the following general rules:

\begin{enumerate}
\item
A minus sign always switches from a LTS to its $c$-dual LTS (indicated by adding a ``prime'' to
the label; if a line contains together with a space its $c$-dual, then the $c$-dual line is omitted). 
\item
For $\K=\C$, a $\C$-linear map $\alpha$ leads to a symmetric space
which is defined over $\C$, whereas a $\C$-antilinear map leads to homotopes of
(pseudo-)Hermitian symmetric spaces. 
\item
Every table contains exactly one line which is of ``groupe type''.
For these cases, the LTS can be directly deduced from the definitions (table in Section
\ref{ClassLieAlgebras}).
\item
In all square matrix cases, there is one symmetric space of the form
``general linear/unitary (resp.\ symplectic)'', arising from the simpler situation of an
algebra with a single involution (Lemma \ref{LemmaB}).
\end{enumerate}

\nin
For reasons of place, we will not spell out computations
for all cases. Let us just
look at some examples: 
in case 1.b of rectangular matrices we
 calculate the LTS corresponding to Line 2 of Theorem \ref{ProjHomotopeTh}:
the symmetric space $\OO(\A,A,\tau) / \OO(\A,A,\tau) \cap \OO(\A,A;\tilde \tau)$ 
where $A$ is the diagonal matrix with blocks
$B$ and $C$,
has LTS
$\A^{(-1,1)}$ with triple product $[X,Y,Z]_A$. Now
$$
\Big[
\begin{pmatrix} 0 & X \cr X^t & 0 \end{pmatrix}, 
\begin{pmatrix} 0 & Y \cr Y^t & 0 \end{pmatrix}, 
\begin{pmatrix} 0 & Z \cr Z^t & 0 \end{pmatrix} \Big]_A =
\begin{pmatrix} 0 & U \cr U^t & 0 \end{pmatrix}, 
$$
with $U =XBY^t C Z + Z BY^t CX -  ( YBX^t C Z + Z BX^t CY)$,
leading to the claimed formula for the LTS. 
The computation for the other rectangular case are similar,  it suffices to replace $\tau$ by
the adjoint matrix w.r.t a Hermitian form on the (skew-)fields $\C$ and $\HHH$.

For the case of the
Siegel-space (label 2.A'), the complex conjugation appearing
in the formula for $\alpha$ stems from the imbedding $\Sym(n,\C) \subset M(2n,2n;\R)$ used
in Theorem \ref{SiegelTh}; the similar case 3.A' arises from the same theorem. 
On the other hand, the imbedding $\Sym(n,\C) \subset M(n,n;\HHH)$ from Theorem
\ref{QuaternionicTh1} leads to the $c$-duals of these two cases (labels 2.A and 3.A).

Under the duality of columns and lines mentioned in the preceding chapter,
the Siegel-space and its analog correspond to the two families 1.1.b' and 1.1.c 
having $V=\Herm(n,\C)$ as underlying space (Theorem \ref{SiegelTh}); similarly
for the two families 1.1.b and 1.1.c' using Theorem \ref{QuaternionicTh1}. 

Finally, tables 3.1 and 2.2 are covered by Theorem \ref{QuaternionicTh2};
mind again the isomorphism (\ref{HermQuat})
 $\Herm(n,\widetilde \HHH) = j \Aherm(n,\HHH)$
(in the complex case we have, of course, $\Aherm(n,\C)=i \Herm(n,\C)$, which explains
why table 1.1 corresponding to this case is longer than the others).
 \EPf
 
 \subsection{Polarized spaces}\label{Sec:Polarized}

 For the following result, recall that a {\em twisted polarized symmetric space} is a symmetric
 space having a local, but not global direct product structure.
 Formally, such spaces have properties very similar to
 (pseudo-) Hermitian symmetric spaces, with complex structures ($J^2 = -1$) replaced by
 polarizations ($J^2 = 1$), see \cite{Be00}.
 If both eigenspaces  of $J$ have equal dimension, one speaks of {\em para-complex
structures} and {\em para-Hermitian symmetric spaces}.
Among the latter are the so-called {\em Cayley type symmetric spaces} (which correspond
precisely to the case when $V^+=V^-=V$ is a {\em Euclidean Jordan algebra}).
Such spaces always have homotopes that are global direct products; we will not list
them here, and focus only on homotopes that are again twisted polarized.

\begin{thm} {\rm (Homotopes of para-Hermitian symmetric spaces)}
The following table contains homotopes of classical para-Hermitian symmetric spaces
that are not globally direct products. The underlying space is in all cases
$\q := V^+ \times V^-$, with Lie triple product
$$
[(X,X'),(Y,Y'),(Z,Z')]_\alpha  = T((X,X'), \alpha(Y,Y'),(Z,Z')) - T((Y,Y'),\alpha(X,X'),(Z,Z'))
$$
with $\alpha:V^+ \times V^- \to V^+ \times V^-$ as below and
$$
T((X,X'),(Y,Y'),(Z,Z')) := (XY'Z + ZY'X,X'YZ' + Z'YX') .
$$
\end{thm}

\msk
\noindent
\begin{tabular}{llll}
label & symmetric space $G_\alpha/H_\alpha$ & $\alpha:V \to V$  &  parameter set 
 \cr
\hline
1.a & 
$\Gl_{2p,2q}( \begin{pmatrix} A & 0 \cr 0 & B \end{pmatrix};\K) /
\Gl_{p,q}(A;\K) \times \Gl_{p,q}(B;\K)$ 
&
$\alpha(X,Y)=(AX^tB,A^tY^tB^t)$
&
$A,B \in M(p,q;\K)$
\cr
1.b & 
$\Gl_{p+q}( \begin{pmatrix} A & 0 \cr 0 & B \end{pmatrix};\K) /
\Gl_p(A;\K) \times \Gl_q(B;\K)$ 
&
$\alpha(X,Y)=(AXA,BYB)$
&
$\begin{matrix} A \in M(p,p;\K), \cr B\in M(q,q;\K) \end{matrix}$
\cr
2. &
$\Sp_n (\begin{pmatrix} 0 & A \cr -A^t & 0 \end{pmatrix};\K) /
\Gl_n(A;\K)$
& $\alpha(X,Y)=(AXA,AYA)$
& $A \in M(n,n;\K)$
\cr
3. &
$\OO_{2n} (\begin{pmatrix} 0 & A \cr A^t & 0 \end{pmatrix};\K) /
\Gl_n(A;\K)$
& $\alpha(X,Y)=(AXA,AYA)$
& $A \in M(n,n;\K)$
\cr
1.1. &
$\UU_{2n} (\begin{pmatrix} 0 & A \cr \overline A^t & 0 \end{pmatrix};\C) /
\Gl_n(A;\C)$
& $\alpha(X,Y)=(AXA,AYA)$
& $A \in M(n,n;\C)$
\cr
3.1. &
$\UU_{2n} (\begin{pmatrix} 0 & A \cr \overline A^t & 0 \end{pmatrix};\HHH) /
\Gl_n(A;\HHH)$
& $\alpha(X,Y)=(AXA,AYA)$
& $A \in M(n,n;\HHH)$
\cr
2.2. &
$\UU_{2n} (\begin{pmatrix} 0 & A \cr \widetilde  A^t & 0 \end{pmatrix};
\widetilde \HHH) / \Gl_n(A;\HHH)$
& $\alpha(X,Y)=(AXA,AYA)$
& $A \in M(n,n;\HHH)$
\end{tabular}

%

\msk
\Pf The proof is a special case of the following more general result:

\begin{thm} {\rm (Homotopes of twisted polarized symmetric spaces)}
The following table contains homotopes of classical twisted polarized symmetric spaces
that are not globally direct products. Using notation as above,
\end{thm}

\msk
\noindent
\begin{tabular}{llll}
label & symmetric space $G_\alpha/H_\alpha$ & $(V^+,V^-)$  &  parameter set 
 \cr
\hline
1. & 
$\Gl_{r+s,r'+s'}( \begin{pmatrix} A & 0 \cr 0 & B \end{pmatrix};\K) /
\Gl_{r,r'}(A;\K) \times \Gl_{s,s'}(B;\K)$ 
&
$(M(r,s'),M(r',s))$
&
$\begin{matrix} A \in M(r',r;\K), \cr B\in M(s',s;\K) \end{matrix}$
\cr
2. &
$\Sp_n (\begin{pmatrix} 0 & A \cr -A^t & 0 \end{pmatrix};\K) /
\Gl_{p,q}(A;\K)$
& $(\Sym(p,\K),\Sym(q,\K))$
& $A \in M(q,p);\K)$
\cr
3. &
$\OO_{p+q} (\begin{pmatrix} 0 & A \cr A^t & 0 \end{pmatrix};\K) /
\Gl_{p,q}(A;\K)$
& $(\Asym(p;\K),\Asym(q;\K))$
& $A \in M(q,p;\K)$
\cr
1.1. &
$\UU_{p+q} (\begin{pmatrix} 0 & A \cr \overline A^t & 0 \end{pmatrix};\C) /
\Gl_{p,q}(A;\C)$
& $(\Herm(p,\C),\Herm(q,\C))$
& $A \in M(p,q;\C)$
\cr
3.1. &
$\UU_{p+q} (\begin{pmatrix} 0 & A \cr \overline A^t & 0 \end{pmatrix};\HHH) /
\Gl_{p,q}(A;\HHH)$
& $(\Herm(p,\HHH),\Herm(q,\HHH))$
& $A \in M(q,p;\HHH)$
\cr
2.2. &
$\UU_{p+q} (\begin{pmatrix} 0 & A \cr \widetilde  A^t & 0 \end{pmatrix};
\widetilde \HHH) / \Gl_{p,q}(A;\HHH)$
& $(\Herm(p,\tilde \HHH),\Herm(q,\tilde \HHH))$
& $A \in M(q,p;\HHH)$
\end{tabular}

\msk
\Pf
Cases 2 and 3 arise from Theorem \ref{ProjHomotopeTh}, and
1.1, 3.1 and 2.2 are treated similarly. Proof for case 1:
Let $\tau:M(p,q;\K) \to M(p,q;\K)$, $X \mapsto I_{r,r'} X I_{s,s'}$.
Clearly $\tau^2 = \id$, and hence, if $\tau(A)=A$ or $\tau(A)=-A$,
$\tau$ is an automorphism of the Lie triple product $[X,Y,Z]_A$, and hence both
eigenspaces are sub-LTS.
Now, both eigenspaces are (as vector spaces) direct products of two spaces of rectangular
matrices having LTS as given in the claim.
\EPf

\begin{rmk} From a Jordan theoretic point of view the latter spaces are all polarized spaces
associated to Jordan pairs  of the form
$(V^+,V^-)=(I,J)$ where $(I,J)$ is a pair of inner ideals in a (simple) Jordan pair $(W^+,W^-)$
(cf.\ Part II).
If $I$ and $J$ are not isomorphic as vector spaces, then the corresponding spaces
are polarized but not para-complex, and they are never semi-simple.
Nevertheless, they have deformations to semi-simple spaces (namely to direct products
of the simple spaces associated to each factor, i.e., to the Jordan pair
$(I,I) \times (J,J)$).  
\end{rmk}

\subsection{Isomorphism classes}

As pointed out in Lemma \ref{strLemma} and the following remark, it is not 
difficult to deduce from the preceding results a classification up to isomorphism.
Generally speaking, isomorphism classes will be parametrized by singling out,
for the matrices $A,B,C$ appearing in the parameter spaces, certain matrices in
normal form: 
 matrices with coefficients $a_{ii}\in \{ 0,1 \}$, $a_{ij}=0$ else,  if $A$ is assumed to be
 rectangular;  
diagonal matrices with coefficients $0,1,-1$ if $A,B,C$ are supposed to be
symmetric or Hermitian; standard skew matrices if $A,B,C$ are assumed to be skew.
Note that the {\em non-degenerate} choices are exactly those corresponding
to the tables from \cite{Be00}, Chapter XII; since these lists are already quite long,
we will not go into details. 

\ssk
Note also that in {\em low dimensional cases} several isomorphisms of spaces occur: in Section
3.7 of Part II we give an overview over such isomorphisms.


\begin{thebibliography}{99}
%




\bibitem[Berger57]{Berger57}
Berger, M, ``Les espaces sym\'etriques non-compactes'', Ann. Ec.\ Norm.\ Sup.\ {\bf 74} (2), 
p. 85--177 

\bibitem[Be00]{Be00}
Bertram, W., {\it The geometry of Jordan- and Lie structures},
Springer Lecture Notes in Mathematics {\bf 1754},  Berlin 2000




\bibitem[Be08]{Be08}
Bertram, W., ``Homotopes and conformal deformations of symmetric spaces.''
 J.\ Lie Theory {\bf 18} (2008), 301 -- 333
 (arxiv: math.RA/0606449)

\bibitem[BeBi]{BeBi}
Bertram, W.\ and Bieliavsky, P.,
``Homotopes of symmetric spaces. II: Structure variety and classification''
(Part II of the present work, arxiv, november 2011)

\bibitem[BeKi09]{BeKi09}{BeKi09b}
Bertram, W.\ and M.\ Kinyon,
``Associative Geometries. I: Torsors, Linear Relations and Grassmannians.
II: Involutions, the Classical Torsors, and their
Homotopes'' , J.\ Lie Theory {\bf 20} (2) (2010), 215--252 and 253--282 
      (arXiv : math.RA/0903.5441 and  math.RA/0909.4438)
      
\bibitem[BDS09]{BDS}
Bieliavsky, P., S.\ Detournay and Ph.\ Spindel,
``The Deformation Quantization of the Hyperbolic Plane'', Com.\ Math.\ Phys.\ {\bf 289} (2),
(2009), 529 -- 559      





\bibitem[Lo69]{Lo69}
  Loos, O., {\it Symmetric Spaces I}, Benjamin, New York, 1969
  



\bibitem[Ma73]{Ma73}
Makarevic, B.O., Open symmetric orbits of reductive groups in symmetric $R$-spaces.
Math.\ USSR Sbornik {\bf 20} (1973), 406 -- 418

\bibitem[Ma79]{Ma79}
Makarevic, B.O., Jordan algebras and orbits in symmetric $R$-spaces.
Trans.\ Moscow.\ Math.\ Soc.\ {\bf 39} (1979), 169 -- 192



\end{thebibliography}
\end{document}